\documentclass[11pt]{amsart}

\newtheorem{thm}{Theorem}[section]
\newtheorem{lem}[thm]{Lemma}

\newtheorem{prop}[thm]{Proposition}

\theoremstyle{plain}
\newtheorem{main}{Theorem}
\newtheorem{corm}[main]{Corollary}
\newtheorem*{blaschke}{Blaschke Conjecture}

\theoremstyle{definition}

\newtheorem*{defn}{Definition}

\numberwithin{equation}{section}

\usepackage{amssymb,latexsym,amscd,amsmath,amsthm}
\usepackage{euscript}
\usepackage{epsfig}

\newcommand{\R}{\mathbf{R}}

\newcommand{\dc}{{\dot{c}}}
\newcommand{\cp}{\mathbf{CP}}
\newcommand{\hp}{{\mathbf{HP}}}

\newcommand{\sph}{\mathbf{S}}

\newcommand{\ml}{\langle}\newcommand{\mr}{\rangle}

\renewcommand{\lim}[1]{\mathop{\underset{#1} {\underset \longleftarrow
{\text{\rm lim}}}}}

\newcommand{\bS}{\bar{S}}

\newcommand{\eps}{\varepsilon}

\newcommand{\sect}{{\bf sec}}
\newcommand{\inj}{\operatorname{inj}}
\newcommand{\diam}{\operatorname{diam}}
\newcommand{\con}{\operatorname{conj}}

\begin{document}

\newcommand{\spacing}[1]{\renewcommand{\baselinestretch}{#1}\large\normalsize}
\spacing{1.14}

\author[K.\ Shankar, R.\ Spatzier, B.\ Wilking]{Krishnan Shankar$^\ast$, Ralf
Spatzier$^{\ast \ast}$, Burkhard Wilking$^{\ast \ast \ast}$}

\title[Spherical rank rigidity and Blaschke manifolds]{Spherical Rank
Rigidity and Blaschke Manifolds}

\thanks{$^\ast$ Supported in part by NSF grant DMS-0103993}
\thanks{$^{\ast \ast}$ Supported in part by NSF grant DMS-0203735}
\thanks{$^{\ast \ast \ast}$ the author is an Alfred P. Sloan Fellow
and was supported in part by NSF}
\subjclass[2000]{53C20, 53C30}
\keywords{Blaschke manifolds, rank one symmetric spaces}

\address{Department of Mathematics, University of Oklahoma, 601 Elm
Ave., Norman, OK 73019.}
\email{shankar@math.ou.edu}

\address{Department of Mathematics, University of Michigan, Ann Arbor,
MI 48109.}
\email{spatzier@umich.edu}

\address{Mathematisches Institut der Uni M\"unster, Einsteinstr. 62,
48149 M\"unster, Germany.}
\email{wilking@math.uni-muenster.de}

\maketitle


\normalsize
\thispagestyle{empty}

\section*{Introduction}

In this paper we define a notion of rank for closed manifolds with
positive upper curvature bound and prove a rigidity result for the
same. More precisely, consider the following definition. If not
explicitly stated otherwise all geodesics are assumed to be
parameterized by arc length.

\begin{defn} Let $M$ be a complete Riemannian manifold with sectional
curvature bounded above by 1.  We say that $M$ has {\em positive
spherical rank} if every geodesic $\gamma\colon [0,\pi]\rightarrow
M$ has a conjugate point at $t=\pi$.
\end{defn}

By the Rauch comparison theorem we know that along any geodesic there
cannot be a conjugate point before $\pi$. The well known equality
discussion implies that for any normal geodesic $c\colon
[0,\pi]\rightarrow M$ there exists a spherical Jacobi field i.e., a
Jacobi field of the form $J(t) = \sin(t) E(t)$ where $E$ is a parallel
vector field (see for instance \cite[Theorem 2.15]{Chavel}).  This
latter characterization is analogous to the notions of
(upper) Euclidean rank and (upper) hyperbolic rank studied by several
people; see below for a more detailed description. In this paper
curvature refers to sectional curvature and is denoted by \sect. The
following is the main result of the paper.

\begin{main}\label{thm1} Let $M^n$ be a complete, simply connected Riemannian
manifold with $\sect \leq 1$ and positive spherical rank. Then $M$ is
isometric to a compact, rank one symmetric space i.e., $M$ is
isometric to  $\sph^n$, $\cp^{\frac{n}{2}}$, $\hp^{\frac{n}{4}}$
or ${\rm Ca}\mathbf{P}^2$.
\end{main}

Note that the condition of $\sect \leq 1$ is not really an
obstruction; any manifold, and in particular any compact manifold,
admits a metric with upper curvature bound 1. So any theorem in this
class must necessarily include an additional assumption on the
geometry of the manifold. There are few general theorems about
manifolds with $\sect \leq 1$; the main theorem and Toponogov's
theorem mentioned in Section 2 are two theorems for such manifolds. We
do not know of any others.

Several notions of `rank' have been studied for manifolds under
suitable curvature assumptions. Historically rank was defined for
symmetric spaces and referred to the dimension of an embedded flat
torus. This evidently descends from the definition of rank for Lie
groups. In our paper we study a more recent notion of rank (also
called the \textit{geometric rank}) first defined in \cite{BBE}
for non-positively curved manifolds. According to \cite{BBE} a
complete Riemannian manifold with $\sect \leq 0$ has
\textit{higher (Euclidean) rank} if along every geodesic $\gamma$
there exists at least one parallel Jacobi field orthogonal to
$\gamma'$. It follows that the 2-plane spanned by this Jacobi
field and $\gamma'$ is extremal.
The following theorem was proved by W.\ Ballmann
(\cite{Ballmann}), and using completely different methods by K.\ Burns
and R.\ Spatzier (\cite{BS}), building on previous work in \cite{BBE}
and \cite{BBS}: Let $M^n$ be a non-positively curved complete manifold
of finite volume. Suppose along every geodesic there exists at least
one parallel Jacobi field. Then the universal cover of $M$ is either a
symmetric space or isometric to a Riemannian product.

The next result was for compact manifolds with $\sect \leq -1$ due to
U.\ Hamenst\"adt; she used a weaker notion of hyperbolic rank by only
assuming that along every geodesic there exists a Jacobi field $J$
such that $\sect(J,\gamma') = -1$ i.e., $J$ and $\gamma'$ span an
extremal curvature 2-plane. She proved the following theorem (cf.\
\cite{Hamenstadt}): Let $M^n$ be a compact manifold with upper
curvature bound $-1$ and hyperbolic rank at least 1. Then $M$ is
isometric to a locally symmetric space. We will refer to this notion
of rank as \textit{(upper) hyperbolic rank}.

In order to exhibit the analogy of these results to the main theorem
we restate it in a slightly weaker form.

\begin{corm} Let $M^n$ be a complete, simply connected Riemannian
manifold with $\sect \leq 1$. Suppose that along every geodesic
$\gamma$, there exists a normal parallel vector field $E$ such that
$\sect(E,\gamma') = 1$. Then $M^n$ is isometric to a compact, rank one
symmetric space.
\end{corm}

Indeed the corollary is an immediate consequence of Theorem~\ref{thm1}
as $\sin(t)E(t)$ is then a Jacobi field along $\gamma$ and
consequently $M$ has positive spherical rank.

Several questions remain open. For instance, one could turn the above
definitions around for (closed) manifolds with suitable lower
curvature bounds and ask whether any rigidity is possible. This is
known to be false if the lower bound is zero but analogous questions
for $\sect \geq -1$ and $\sect \geq 1$ remain untouched. We refer the
reader to Table 1 where some of the known results are presented; the
table is not meant to be a survey rather a point of departure for
further investigations. In this paper we only deal with the case of
spherical rank for upper curvature bound 1.

The paper is organized into three sections.  In Section 1 we show that
positive spherical rank implies that the manifold is a so-called
\textit{Blaschke manifold}. A Blaschke manifold is a Riemannian
manifold with the property that its injectivity radius equals its
diameter. Note that the definition has no curvature assumptions. For
an excellent and rather complete treatment of Blaschke manifolds see
\cite{Besse}. The study of these manifolds has a rich history
motivated by the following open problem.

\begin{blaschke} Let $M$ be a Riemannian manifold such that $\inj M =
\diam M$. Then $M$ is isometric to a compact, rank one symmetric space
(CROSS).
\end{blaschke}

Once we have established that $M$ is Blaschke the remaining step may
be regarded as a special case of the Blaschke conjecture. In Section 2
we show that a Blaschke manifold with upper curvature bound 1 and
injectivity radius $\pi$ must be isometric to a CROSS. This latter
result has already been proved by V.\ Rovenskii and V.\ Toponogov in
\cite{RT}; they prove this using comparison arguments. For the sake of
completeness we present a shorter proof whose arguments may be useful
in other contexts. In the final section of the paper we give an
example, the Berger spheres, which shows that the conclusion of the
main theorem fails for a weaker notion of spherical rank, namely the
analogue of Hamenst\"adt's notion of rank.  More precisely, we show that
there are non-symmetric, simply connected, compact Riemannian
manifolds with upper (respectively lower) curvature bound 1, such that
along every geodesic $\gamma$ there exists a normal Jacobi field $J$
with  $\sect(\gamma'(t),J(t))=1$ for all $t$.
\medskip

It is a pleasure to thank Karsten Grove and Wolfgang Ziller for
several helpful discussions.

\section{Positive Spherical Rank implies Blaschke}

Let $M$ be a complete, simply connected Riemannian manifold with
$\sect \leq 1$ and positive spherical rank. By assumption, every
geodesic hits its first conjugate point at $\pi$ and therefore the
diameter of $M$ is bounded above by $\pi$. In order to show that $M$ is
a Blaschke manifold we only have to verify that the injectivity radius
is at least $\pi$ since we always have $\inj M \leq \diam M$. Since
the conjugate radius of the manifold is $\pi$, it suffices to show
$\inj M = \con M$.

Consider the special case where $M$ is even dimensional and
positively curved. Then by Klingenberg's injectivity radius
estimate $\inj M = \pi$, and hence $M$ is Blaschke. This observation,
in fact, was the beginning for our investigations.

We now outline the argument for the general case.  The starting point
is the well-known generalization \cite[Lemma 5.6]{CE} of an
injectivity radius estimate of Klingenberg \cite[Lemma 4]{Klingenberg}
that for $M$ compact, $\inj M$ is the smaller of $\con M$ and half the
length of a shortest closed geodesic. We will argue by contradiction,
and suppose that the length of some closed geodesic is less than
$2\pi$. Under these assumptions, it follows from Morse theoretic
arguments that there exists a closed geodesic $\gamma$ of length $2
\pi$ and index 1.\footnote{Throughout this paper we always mean index
in the free loop space.}  For the actual argument it is important that
$\gamma$ satisfies a slightly stronger condition; see Lemma~\ref{lem:
deg Morse}.

Moreover, we will show that $\gamma$ is contained in a totally
geodesic, isometrically immersed 2-sphere of constant curvature 1. The
next step is to show that the same is true for all geodesics in the
manifold and hence all geodesics are closed. Then in Section 1.4 we
will show, by applying the index parity theorem \cite{Wilking} that if
all geodesics are closed, then they all have length at least $2
\pi$. This contradiction finishes the proof that $M$ is indeed a
Blaschke manifold.

\subsection{Preliminaries} \hfill

%

First we state a useful generalization of Klingenberg's long homotopy
lemma due to U.\ Abresch and W.\ Meyer \cite{Abresch-Meyer}.

\begin{lem}[Long Homotopy Lemma] Let $M$ be a compact Riemannian
manifold and $c$ a closed curve in $M$ which is the union of at most
two geodesic segments such that $l(c) < 2\con M$. Suppose $c= c_0$ is
homotopic to a point via a continuous family of rectifiable closed
curves $c_t,\, 0 \leq t \leq 1$. Then some $c_s$ has length $l(c_s)
\geq 2 \con M$.
\end{lem}

Now we adapt the second comparison theorem of Rauch to get the
next proposition.  We will need the following lemma in the proof.

\begin{lem}Let $M$ be a complete manifold with $\sect \leq 1$, and
let $\sph ^2$ denote the 2-sphere of constant curvature 1. Suppose $X$
is a normal Jacobi field along a geodesic $b$ of length at most $\pi$ such
that $\|X(0)\|=1$ and $\ml X',X\mr =0$. Let $Y$ be a normal
Jacobi field along a
geodesic in $\sph ^2$ such that $\|Y(0)\|=1$ and $Y' (0) =0$.

 Then
$\|X(s)\| \geq \|Y(s)\|$ for $0 \leq s \leq \pi/2$.  Moreover, if
$\|X(s_1)\| = \|Y(s_1)\|$ for some $0 < s_1 \leq \pi/2$, then $\sect
(X(s), b' (s))=1$ and $\|X(s)\|=\|Y(s)\|$ for all $0 \leq s \leq s_1$ .
\end{lem}

\begin{proof} Since the upper curvature bound is 1, we see from
Cauchy-Schwarz and straightforward differentiation that
\[ \| X \| '' + \| X \| \geq 0 .\] Let $a(s)$ be the function such
that $ \| X \| '' + a(s) \| X \| = 0$. Note that $a(s) \leq \sect
(X(s),b' (s)) \leq 1$.  By the Sturm comparison theorem, it follows
that $\| X (s) \| \geq \| Y (s) \|$ on the closed interval
$[0,\pi]$. We refer to \cite[p.\ 238]{docarmo} for a Sturm comparison
theorem with different initial conditions. The same proof however
applies equally well in our situation.  Moreover if $\| X (s_1) \| =
\| Y (s_1) \|$, then $a(s) =1$ for all $s \leq s_1$. In particular,

$\sect (X(s),b'(s)) =1$ for $0 \leq s \leq s_1$.
\end{proof}

An isometrically immersed surface (with piecewise smooth geodesic
boundary) of constant curvature 1 will be called a \textit{spherical
slice}.
\medskip

Let $M$ be a complete, simply connected manifold with $\sect \leq
1$. Suppose $c$ is a geodesic of length $\pi$ between two points $p$

and $q$ on $M$ and suppose $c_s$ is a smooth variation of $c$ by
curves connecting $p$ and $q$ such that $l(c_s) \leq \pi$ for all
$s$. We pick a curve $c_\sigma$ close to $c$ and assume that
$c_\sigma$ is not a reparameterization of $c$.

%

\begin{prop}\label{prop: sph slice} Given $c$ and $c_\sigma$ as above,
they span a totally geodesic spherical slice.
\end{prop}

\begin{proof} If we choose $\sigma$ sufficiently small
we can find a normal vector field $Z(t)$ along $c$ such that
$c_\sigma(t) = \exp(Z(t))$ after possibly reparameterizing
$c_{\sigma}$. We may assume $\|Z(t)\|<\pi/2$.  This yields a proper
variation $f(s,t) = \exp(s\cdot Z(t))$ of $c$. By construction, the
curves $f_t(s)$ ($t$ fixed) are geodesics which in turn implies that
the vector field $X_t(s) = \frac{\partial f}{\partial t}$ is a Jacobi
field along $f_t(s)$.  It should be clear that $f_t$ is not
necessarily parameterized by arc length. However, according to our
conventions $c$ is and thus $\|X_t(0)\| = \|c'(t)\|=1$.

Let $\tilde{c}$ be a geodesic of length $\pi$ on $\sph^2$, the
2-sphere with constant curvature 1. Consider the following variation
of $\tilde{c}$,
$$
    g(s,t) = \exp\bigl(s\cdot \|Z(t)\|\cdot \widetilde{E}(t)\bigr),
$$
where $\widetilde{E}$ is a unit parallel field along $\tilde{c}$
orthogonal to $\tilde{c}'$. Notice that $g_s(t)=g(s,t)$
is a proper variation as well.
 Let $Y_t(s) = \frac{\partial g}{\partial t}$ be the
vector field along the geodesics $g_t(s)$. Then
$$
\|Y_t(0)\| = \biggl\|\frac{\partial g}{\partial t}(0) \biggr\| =
\|\tilde{c}'(t)\| = 1.
$$

We would like to apply the Sturm comparison theorem to the vector
fields $X_t(s)$ and $Y_t(s)$. To do this we need to estimate the
derivatives of these vector fields. A straightforward calculation
yields
\begin{align*}
\|Y_t(0)\|' &= \frac{\partial}{\partial s}(\|Y_t(s)\|)_{s=0} = 0 , \\
\|X_t(0)\|' &= \frac{\partial}{\partial s}(\|X_t(s)\|)_{s=0} = 0.
\end{align*}
As always each of the Jacobi fields $Y_t$ and $X_t$ can be decomposed
into the sum of a normal and a tangential Jacobi field along $f_t$
resp. $g_t$. The tangential parts of the Jacobi fields are given by
$m_t\, s\, f_t'(s)$ respectively by $m_t\, s\, g_t'(s)$ with
$m_t=\tfrac{\partial}{\partial t}\log\bigl(\|Z(t)\|\bigr)$.  Since the
geodesics $g_t$ and $f_t$ have the same speed it follows from the
previous lemma that the norm of the normal part of $X_t$ is bounded
above by the norm of the normal part of $Y_t$.  Combining the two
statements we get $\|X_t(s)\| \geq \|Y_t(s)\|$ and so
 \[
\pi \leq \int \|Y_t(s)\|\, dt \leq \int \|X_t(s)\|\, dt
\] for
all $s$. By construction equality holds at $s=1$.  Notice that this
implies in particular that $c_{\sigma}$ is a geodesic up to
parameterization, as otherwise one could have replaced $c_{\sigma}$ by

a nearby curve of length $<\pi$.  The equality discussion implies
$\|X_t(s)\|=\|Y_t(s)\|$ for $s\in [0,1]$. This shows that the strip
parameterized by $f$ is intrinsically isometric to the strip
parameterized by $g$. Furthermore the equality discussion also shows
that the ambient curvature of the slice defined by $f$ is $1$ as well.
By the Gauss Lemma we have the basic relation
$$
0=\sect_{\,\text{intrinsic}} - \sect_{\,\text{ambient}} = \langle
B(X_1,X_1), B(X_2,X_2) \rangle - \| B(X_1,X_2)\|^2,
$$
where $X_1 = \frac{\partial f}{\partial t}$ and $X_2 = \frac{\partial
f}{\partial s}$ are linear independent vector fields on the slice. By
construction, the curves $f_t(s)$ ($t$ fixed) are geodesics so
$B(X_2,X_2) = 0$. Therefore, $B(X_1,X_2)$ vanishes as well. It remains
to show that $B(X_1,X_1)=0$.  Notice that in the 'model' slice
parameterized by $g$ the geodesics of length $\pi$ connecting the end
points of $\tilde{c}$ pass through every point of the slice.  Since
the two slices are intrinsically isometric the same holds for the
slice parameterized by $f$.  Notice that these intrinsic geodesic have
to be geodesics of the ambient manifold as well because otherwise we
could find in $M$ nearby curves which are strictly shorter.
\end{proof}

\subsection{Existence of a closed geodesic of length at least $2\pi$
and index 1} \hfill

The next lemma is a consequence of the long homotopy lemma; the first
part is well known.

\begin{lem}\label{lem: deg Morse} Let $M$ be a complete, simply connected
compact manifold with $\sect \leq 1$ and injectivity radius less than $\pi$.
Then there is a closed geodesic $c\colon [0,\ell]\rightarrow M$ of
length $\ell \ge 2 \con M \ge 2\pi$ whose index in the free loop space
of $M$ is at most 1.

Furthermore there is no free homotopy $c_s(t)$ with $s\in [0,1]$ and
$t\in[0,\ell]$ such that each of the following statements is true.

\begin{enumerate}
\item[(\textit{i})] $c_s$ is a closed geodesic of length $2\con M$ and
$c_0=c$.

\item[(\textit{ii})] The index of the closed geodesic $c_s$ in the
free loop space of $M$ is at least $1$, $s\in [0,1]$.

\item[(\textit{iii})] The index of the closed geodesic $c_1$ is at
least $2$.
\end{enumerate}
\end{lem}

For the proof of the above lemma we will apply the standard degenerate
Morse Lemma, see for example \cite{GM69}.

\begin{lem}\label{morse lem} Let $B$ be a manifold of dimension $b$,
$E\colon B\rightarrow \R$ a smooth proper function
and let $p\in B$ be a critical point of $E$.
Then we can find a neighborhood $U$ of $p$
and a map $x\colon U\rightarrow V\subset \R^b$ with $x(p)=0$ such
that
\[
E = E(p) -x_1^2-\ldots -x_{\lambda}^2+ x_{\lambda+1}^2+\ldots +x_{b-d}^2
+h(x_{b-d+1},\ldots, x_{b})
\]
where $\lambda$ denotes the index of $p$,  $d$ the nullity of $p$
and $h$ is a smooth function.
\end{lem}

Notice that any critical point of $E$ in $U$ is necessarily contained
in $L:=x^{-1}(0\times \R^d)$. After replacing $U$ by a smaller
neighborhood we may assume that $V$ is a bounded convex set.  For the
proof of Lemma~\ref{lem: deg Morse} we make the following observation:
Suppose that $\lambda>0$.  Let $p^i$ be a sequence of points
converging to $p$ with $E(p^i)<E(p)$, $h(t)\in U$ a path of critical
points with $E(h(t))=E(p)$ and $h(0)=p$, $t\in [0,1]$. Then there is a
path $h_i(t)$ with $p^i=h_i(0)$, $E(h_i(t))<E(p)$ such that $h_i(1)$
converges to $h(1)$.  In order to construct $h_i$ we will identify $U$
with $V$ via $x$.  Consequently we write $p^i_j$ instead of
$x_j(p^i)$.  First take a path given by the straight segment from
$p^i$ to another point $q^i$ with $q^i_j=p^i_j$ for $j\ge 2$ and
$|q^i_1|> \eps$, where $\eps>0$ is a number which we can chose
independent of $i$.  Next consider the path from $q^i$ to
$\tilde{q}_i:=(q^i_1,0,\ldots,0)$ given by a straight segment.  Since
$p_i$ converges to $0$ it is easy to see that the energy $E$ along
this path stays strictly below $E(p)$ for almost all $i$.  Next
consider the path $\tilde{h}_i(t)= \tilde{q}^i + h(t)$ from
$\tilde{q}^i$ to $\tilde{q}^i + h(1)$ along which the energy $E$ is
constant.  Finally, along the straight line from $h(1)+\tilde{q}^i $
to $h(1)+\tfrac{\tilde{q}^i}{i}$ the energy stays strictly below
$E(p)$. Thus we may chose $h_i$ as the composition of these paths for
almost all $i$. Finally we can define $h_i$ as the point curve for the
finitely many remaining $i$.

\begin{proof}[Proof of Lemma~\ref{lem: deg Morse}.] In this
proof all curves are parameterized on $[0,1]$.  As usual we consider
the energy functional on the free loop space of $\Omega M$ of $M$
i.e., we define the energy of a piecewise smooth loop $c\colon
[0,1]\rightarrow M$ as
\[
E(c):=\tfrac{1}{2}\int_0^1\|\dot{c}(t)\|^2\,dt.
\]
For any value of $e\in (0,\infty)$ we let $\Omega M^{<e}$
(respectively $\Omega M^{\le e}$ ) denote the loops in $\Omega M$ of
energy $<e$ (respectively $\le e$).

By Klingenberg's general injectivity radius estimate there is a closed
geodesic of length $2\inj M <2\pi$. Furthermore the long homotopy
lemma tells us that there is no free null homotopy of this geodesic
contained in $\Omega M^{<e_0}$ with $e_0=2(\con M)^2$.  Consequently
$\Omega M^{<e_0}$ has at least two connected components.

We now assume, on the contrary, that the statement of the Lemma is
false.  The first step is to verify that $\Omega M^{\le e_0}$ is
connected.  In fact if $\Omega M^{\le e_0}$ were not connected, then
we could find an $\eps>0$ such that $\Omega M^{< e_0+\eps}$ is not
connected either.  Since the statement of the lemma is assumed to be
false, any closed geodesic of energy $>e_0$ has index at least $2$.
Thus it follows by the usual degenerate Morse theory argument (namely
approximating $E$ by Morse functions) that the free loop space itself
is not connected either.  This is a contradiction as $M$ is simply
connected. Hence, $\Omega M^{\le e_0}$ is connected.

As usual given an $e_1$ one can find partition $ 0<t_1<\cdots <t_k<1$
of the unit interval such that for all $e\le e_1$ the sub level $\Omega
M^{\le e}$ is homotopically equivalent to the subset of broken
geodesics $B^{\le e}$ contained in $\Omega M^{\le e}$, whose points of
non differentiability are points in the partition.  We put $e_1=e_0
+1$ and fix a sufficiently fine partition. Then $B^{< e_1}$ is a
finite dimensional submanifold and if we restrict the energy function
to $B^{< e_1}$, then the critical point as well as the indices do not
change.  We have shown that $B^{\le e_0}$ is connected while $B^{<
e_0}$ is not.

Let $C$ denote the set of closed geodesics of length $2(\con M)$ and
put $S= B^{< e_0}\cup C$.  In other words $S$ is obtained from
 $B^{\le e_0}$ by removing all non-critical points from the boundary.  Since
$B^{\le e_0}$ is connected it is easy to see that $S$ is connected as
well.

Let $S_1$ be an open and closed subset of $B^{< e_0}$ and suppose that
neither $S_1$ nor its complement $S_2:=B^{< e_0}\setminus S_1$ is
void.  Let $\bS_i$ denote the closure of $S_i$ in $S$.  By
construction $\bS_1\cup\bS_2=S$ and $\bS_1\cap \bS_2$ is a nonempty
subset of $C$. We claim that if $c\in \bS_1\cap \bS_2$ then
$C_0\subset \bS_1\cap \bS_2$ where $C_0$ is the path connected
component of $c$ in $C$.  Since the index of any critical point in $C$
is at least $1$ this is an immediate consequence of the observation
that we made after Lemma~\ref{morse lem}.

Let $C'\subset C$ denote the set of closed geodesics of index $\ge 2$.
By assumption each path connected component can be represented by a
geodesic of index $\ge 2$.  By the previous argument $\bS_1\cap \bS_2$
has a nontrivial intersection with $C'$.  But this shows that $S':=
B^{< e_0}\cup C'$ is connected as well.  Since all points in $C'$ have
index at least 2 this implies as before that $ B^{< e_0}$ is connected
which is a contradiction.
\end{proof}

\subsection{All geodesics in $M$ are closed} \hfill

We will assume from now on that $M$ is a complete, simply connected
manifold with $\sect \leq 1$ and positive spherical rank.  In this
subsection we want to prove that all geodesics of $M$ are
closed. There is nothing to prove if $M$ is Blaschke. Thus we may
assume $\inj M < \pi$. Then there is a geodesic $\gamma$ satisfying
the conclusion of Lemma~\ref{lem: deg Morse}.  Using the existence of
$\gamma$ as a starting point we will show that all geodesics are
contained in a totally geodesic immersed $2$-sphere of constant
curvature one.

Consider the set $\mathcal{R}$ of geodesic segments $c:[0,\pi]
\rightarrow M$ for which $\pi$ is a conjugate point with multiplicity
1. Note that $\mathcal{R}$ is an open set in the set of all geodesic
segments of length $\pi$.  Indeed, the multiplicity of a limit of
geodesics of length $\pi$ can only go up.  Moreover, since the
spherical rank is positive, every such geodesic has multiplicity at
least 1. Finally $\mathcal{R}$ is not empty as it contains
$\gamma_{|[0,\pi]}$.

\begin{lem}\label{lem: sph slice}
Suppose $c \in \mathcal{R}$ and $J$ is a Jacobi field along $c$ that
vanishes at 0 and $\pi$ with $\|J'(0)\|=1$. Put $v:=\dc(0)$ and
$w=J'(0)$.

Then there is a unique maximal number $s_m\in (0,\pi]$ such that
\[
h(s,t)=\exp\bigl(t(\cos(s)v+\sin(s)w)\bigr),
\mbox{ with  $t\in [0,\pi]$ and $s\in [-s_m,s_m]$}
\]
parameterizes a totally geodesic immersed spherical slice of constant
curvature 1.  If $s_m<\pi$ then one of the boundary geodesics $h(\pm
s_m,t)$ is not contained in $\mathcal{R}$.  If $s_m=\pi$ then the
image of $h$ is a totally geodesic immersed $2$-sphere.
\end{lem}

\begin{proof} We first want to show that we can indeed chose $s_m>0$.
Let $v\in T_pM$ be the initial vector of $c$. Denote by $T^1_p M$ the
unit sphere in $T_p M$. Consider the map $\phi: T^1_p M \rightarrow M$
given by $w \mapsto \exp( \pi w)$, and set $q= \phi(v)$.  Since the
spherical rank is positive, $\phi$ has a singular differential
everywhere.  At $v$ the kernel of the differential is precisely one
dimensional since $c\in \mathcal{R}$. It follows that $\phi$ is
of constant rank in a neighborhood of $v$. By the implicit function
theorem, the fibers $\phi ^{-1} (x)$ are 1-dimensional submanifolds
for $x$ in a neighborhood of $q$. The curve $\phi ^{-1} (q)$ defines a
variation of geodesics $f(s,t)$ of $c(t)$ of length $\pi$ with
constant starting and ending point. From Proposition~\ref{prop: sph
slice} we deduce that $c$ is contained in a spherical slice.  Since
there is up to constant factor only one Jacobi field along $c$ that
vanishes at $0$ and $\pi$ this spherical slice is necessarily contains
$h(s,t)$ for all $(t,s)\in [0,\pi]\times [-s_m,s_m]$ provided that
$s_m>0$ is chosen sufficiently small.  This proves $s_m>0$ and clearly
we may choose $s_m\in (0,\pi]$ maximal.

Consider next the case of $s_m<\pi$.  Notice that
$J_{\pm}(t)=\tfrac{\partial h}{\partial s}(\pm s_m,t)$ is a Jacobi
field along the boundary geodesic $c_\pm(t)=h(\pm s_m,t)$. If $c_\pm
\in \mathcal{R}$, then the previous argument shows that we can
increase $s_m$ contradicting our choice of $s_m$.
\end{proof}

Consider again the geodesic $\gamma$ from Lemma~\ref{lem: deg Morse}.
Since $M$ has positive spherical rank and $\gamma$ has index at most
1, it follows that $l(\gamma) = 2 \pi$. We parameterize $\gamma$ on
the interval $[-\pi, \pi]$.

\begin{prop} The closed geodesic $\gamma$ of length $2 \pi$ and index
1 is contained in a totally geodesic, isometrically immersed $\sph^2$
of constant curvature 1. \end{prop}

\begin{proof} We construct two continuous vector fields $X_+$ and
$X_-$ along $\gamma$ that are defined as follows:

(\textit{i}) $X_+(t) = 0$ for all $t<0$ and $X_-(t) = 0$ for all
$t>0$.

(\textit{ii}) $X_+(t)$ is a non-vanishing Jacobi field on $[0, \pi]$
and $X_-(t)$ is a non-vanishing Jacobi field on $[-\pi, 0]$.

Note that the index form of $\gamma$ restricted to the two dimensional
subspace spanned by $X_+$ and $X_-$ is $0$. Since the index of
$\gamma$ is at most 1, it follows that the two dimensional space
spanned by $X_+$ and $X_-$ must contain a Jacobi field $J = aX_+ +
bX_-$. By the equality discussion of the Rauch comparison theorem,
$X_{\pm}$ looks like $\sin(t) E_{\pm}$ on the intervals where they are
non-zero; here $E_{\pm}$ are parallel vector fields. Since $X$ is a
smooth Jacobi field on $[-\pi, \pi]$ it follows by computing $X'$ that
$J(t) = \sin(t) X(t)$, where $X = aE_+ = bE_-$ i.e., $J$ is a periodic
Jacobi field. Moreover, $X(t)$ is a closed, parallel vector field
along $\gamma$ such that $\sect(\gamma'(t), X(t)) = 1$ for all $t$.

The vector fields $\sin(t) X(t)$ and $\cos(t) X(t)$ are periodic
Jacobi fields along $\gamma$. From the previous lemma it follows that
there is an $\epsilon >0$ such that
\[
\exp(s\, X(t)),\mbox{   $t \in [-\pi, \pi]$ and $s\in
[-\epsilon, \epsilon]$}
\]
parameterizes a totally geodesic spherical tube. Notice that there
are lots of closed geodesics in this tube. Every one of the closed
geodesics in the tube is homotopic to $\gamma$ via a homotopy
$\gamma_s$ satisfying the first two conditions of Lemma~\ref{lem: deg
Morse}. By the same lemma it follows that the third condition must be
violated i.e., each of the closed geodesics in the tube must have
index one in the free loop space. Therefore, there is no obstruction
to increase $\eps$.  In other words, we may choose $\eps=\pi/2$ and
thus $\gamma$ is contained in a totally geodesic immersed $\sph^2$ of
constant curvature 1.
\end{proof}

It is important to notice that each of the closed geodesics in the
constant curvature $\sph^2$ constructed above has index 1 in the free
loop space. This implies that along every geodesic $t=\pi$ is a
conjugate point with multiplicity 1.

\begin{prop} Suppose $M$ has $\sect \leq 1$ and positive spherical
rank. Then all geodesics in $M$ are closed.
\end{prop}

\begin{proof} Consider the following subsets of $T^1 M$, the unit
tangent bundle of $M$.
$$
S_1 = \biggl\{ v \in T^1 M \Bigm|
\begin{array}{l} \mbox{ $v$ tangent to a totally geodesic immersed

$\sph^2_v$ and all geodesics}\\\mbox{ $c\colon [0,\pi]\rightarrow \sph^2_v$
have a conjugate point of multiplicity 1 at $\pi$.}
\end{array}\biggr\}
$$
$$
S_2 = \biggl\{ v \in T^1 M \Bigm| \begin{array}{l} \mbox{ $v$
tangent to a totally geodesic immersed $\sph^2_v$ and all
geodesics}\\\mbox{ $c: [0, 2\pi] \rightarrow \sph^2_v$ have index 1 in
the free loop space $\Omega M$.}
\end{array}\biggr\}
$$

Let $v_0 \in T^1M$ denote the initial velocity vector of the closed
geodesic $\gamma$ of length $2\pi$. Then one can see that $S_2 \subset
S_1\subset T^1 M$.  Furthermore $S_2$ is non-empty since it contains
$v_0$ and $S_2$ is closed.

Next we claim that $S_1$ is open.  Let $w\in S_1$ and let $\sph^2_w$
be as in the definition of $S_1$.  Suppose a sequence $w_i\in T^1M$
converges to $w$.  For $i$ sufficiently large there is a unique
spherical Jacobi field $\sin(t)X_i(t)$ along the geodesic
$c_i(t)=\exp(tw_i)$, $t\in [0,\pi]$.  Suppose for a moment that $w_i$

is not tangent to a totally geodesic immersed 2 sphere.  By
Lemma~\ref{lem: sph slice} $w_i$ is tangent to a spherical slice such
that one of the boundary geodesics is not contained in $\mathcal{R}$.
Since a subsequence of the boundary geodesics converges to a geodesic
in $\sph^2_w$ and $\mathcal{R}$ is open this is impossible.  In other
words $w_i$ is tangent to a totally geodesic immersed sphere
$\sph_{w_i}^2$. Since the geodesics in $\sph^2_{w_i}$ converge to
geodesics in $\sph^2_{w}$ we deduce that $w_i\in S_1$ for almost all
$i$.

Next we finish up the proof of the proposition by showing $S_2 = M$.
Suppose, on the contrary that $S_2\neq M$. Choose a path $h(s)\in
T^1M$ with $h(0)=v_0$ and $h(1)\in T^1M\setminus S_2$.  Since $S_1$ is
an open neighborhood of the closed set $S_2$ we may assume that
$h(s)\in S_1$ for all $s\in [0,1]$.  Furthermore we may assume that
$c_1(t)=\exp(th(1))$ ($t\in [0,2 \pi]$) is one of the closed geodesics
in $\sph^2_{h(1)}$ of index at least $2$.  Thus
\[
 c_s(t)=\exp(th(s)) \,\,\mbox{for  $t\in[0,2\pi]$, $s\in [0,1]$}
\]
defines a homotopy of closed geodesics of length $2\pi$ satisfying all
three conditions of Lemma~\ref{lem: deg Morse} which is a
contradiction.
\end{proof}

\subsection{$M$ is a Blaschke manifold} \hfill

We now show that $M$ is Blaschke i.e., $\inj(M)\ge \pi$. By the
previous subsection all geodesics of $M$ are closed.  This enables us
to apply the following index parity theorem (cf.\ \cite{Wilking}):

\begin{thm}[Wilking] Let $M^n$ be an oriented Riemannian manifold all
of whose geodesics are closed, and let $c:[0,1] \rightarrow M$ be a
closed geodesic. Then the index of $c$ in the free loop space of $M$
is even if $M$ is odd-dimensional and it is odd if $M$ is
even-dimensional.
\end{thm}

\begin{proof}[Proof that $M$ is Blaschke] We argue by contradiction
and assume that $\inj M < \pi$. Then by the generalized injectivity
radius estimate of Klingenberg, there exists a shortest closed
geodesic $\alpha$ of length $2\,\inj M$.  By the Long Homotopy Lemma
(Lemma 1.1) we know that $\alpha$ is not freely null homotopic in the
space of all curves of energy less than $2\pi^2$ (or all curves of
length shorter than $2\pi$). In particular, it follows that the curves
of length less than $2\pi$ form a disconnected set such that $\alpha$
and the point curve lie in distinct components. On the component
containing $\alpha$, the energy functional attains a minimum at
$\alpha$ and hence $\alpha$ has index 0 in the free loop
space. However, if $M$ is even-dimensional, then the index of $\alpha$
must be odd by the index parity theorem which leads to a
contradiction.  If $M$ is odd dimensional consider the closed geodesic
$\gamma$ constructed in Section 1.2. By construction, $\gamma$ has
index exactly 1. Once again this contradicts the index parity theorem
as the index of $\gamma$ is required to be even. \end{proof}

\section{A special case of the Blaschke conjecture}

We have shown that if $M$ is a complete, simply connected, Riemannian
manifold with $\sect \leq 1$ and with positive spherical rank, then
$\inj M = \diam M = \pi$ i.e., $M$ is a Blaschke manifold with
extremal diameter (and injectivity radius). In this section we
complete the proof of the main theorem by proving the following
proposition which is special case of the Blaschke conjecture.

\begin{prop} Let $M$ be a simply connected Blaschke manifold with
$\sect \leq 1$ and extremal value of diameter (and injectivity radius)
equal to $\pi$. Then $M$ is isometric to a compact, rank one symmetric
space. \end{prop}

As we noted in the introduction, the above result has already been
proved by Rovenskii and Toponogov in \cite{RT}; they also use
Toponogov's theorem below. The proof given here has the slight virtue
of being shorter. To be more precise we reduce the problem to two
older theorems: one due to V.\ Toponogov (see \cite{Toponogov}) and
the other due to M.\ Berger (see \cite{Berger}).

\begin{thm}[Toponogov] Let $M$ be a complete, simply connected,
Riemannian manifold such that $\sect_M \leq 1$. Suppose $M$ contains a
closed geodesic $\gamma$ of length $2\pi$ and index $k-1$. Then
$\gamma$ is contained in an isometrically embedded, totally geodesic,
sphere $\sph^k$ of constant curvature 1.
\end{thm}

It should be noted that the proof of Toponogov's theorem is not very
hard in the special case that $M$ is Blaschke. In that case the map

$f_p\colon T^1_pM\rightarrow M$, $v\mapsto \exp(\pi v)$ has constant
rank for all $p\in M$, namely the rank equals $n-k$, where $k-1$ is
the index of a closed geodesic of length $2\pi$ or equivalently the
multiplicity of the conjugate point at $\pi$.  Thus the fibers of
$f_p$ are submanifolds and using Proposition~\ref{prop: sph slice} it
is easy to see that the fibers are great spheres of dimension $k-1$.
Furthermore one can use Proposition~\ref{prop: sph slice} to see that
$\exp(\R\cdot f_p^{-1}(q))$ is a totally geodesic sphere of dimension
$k$ and of constant curvature 1.
\medskip

Before we state Berger's theorem some notation is required. An
$SC_{2a}$-manifold is one in which every geodesic is simply closed and
periodic with period $2a$. It is well known that a Blaschke manifold

with $\inj = \diam = a$ is an $SC_{2a}$-manifold (cf.\
\cite[Chapter 7]{Besse}). In our situation, we have a Blaschke
manifold which happens to be an $SC_{2\pi}$-manifold, so every geodesic
is simply closed with period $2\pi$.

Given two points $p,q$ at distance $\pi$ on a Riemannian manifold $M$,
let $\Sigma_\pi(p,q)$ denote the set of all shortest geodesics from
$p$ to $q$. It is shown in \cite{Besse} that in this case $\Sigma_\pi
(p,q)$ is homeomorphic to a sphere $\sph^k$, where $k-1$ is the index
of a closed geodesic through $p$ and $q$. If for all tuples $(p,q)$

with $d(p,q)=\pi$, the set $ \Sigma_\pi(p,q)$ is totally geodesic,
then, following Berger \cite{Berger}, $M$ is called a \textit{totally
geodesic Blaschke manifold}.

\begin{thm}[Berger] Let $M$ be a simply connected, totally geodesic Blaschke
manifold. Then $M$ is isometric to a compact, rank one symmetric space
i.e., isometric to $\sph^n$, $\cp^{\frac{n}{2}}$, $\hp^{\frac{n}{4}}$
or ${\rm Ca}\mathbf{P}^2$. \end{thm}

\begin{proof}[Proof of Proposition 2.1] If a complete, simply connected
Riemannian manifold has $\sect \leq 1$ and positive spherical rank,
then it is an $SC_{2\pi}$ Blaschke manifold. So every geodesic in $M$
is simply closed, has length $2\pi$ and index at least 1.

Pick any geodesic $\gamma$ of length $2\pi$ and index $k-1\geq 1$ and
pick two points $p,q$ on $\gamma$ which are $\pi$ apart. By
Toponogov's theorem $\gamma$ is contained in a totally geodesic,
isometrically embedded $\sph^k$ of constant curvature 1. By
construction we have $\sph^k \subset \Sigma_\pi (p,q)$. But
$\Sigma_\pi (p,q)$ is also a $k$-dimensional sphere because of the
index estimate on $\gamma$. Moreover, $\Sigma_\pi (p,q)$ is connected
which implies $\sph^k = \Sigma_\pi (p,q)$ and $M$ is a totally
geodesic Blaschke manifold. By Berger's theorem $M$ must be isometric
to a CROSS. \end{proof}

\section{Some Examples}

In this section we explore another notion of spherical rank that is
analogous to Hamenst\"adt's notion of hyperbolic rank. More precisely,
consider the following:

\begin{defn} Let $M^n$ be a compact Riemannian
manifold.  Suppose along every geodesic $\gamma(t)$ in $M$ there
exists a normal Jacobi field $J(t)$ such that $\sect(\gamma'(t),J(t))
= 1$.  If $\sect \le 1$, we say that $M^n$ has \textit{weak upper
spherical rank} at least 1. If $\sect \ge 1$, we say that $M^n$ has
\textit{weak lower spherical rank} at least 1.
\end{defn}

In the case of stronger notions of rank we have seen various rigidity
results, most of them implying that the universal cover must be
locally isometric to a symmetric space. We direct the reader to Table
1 for some of the known results.  As is indicated there (metric)
rigidity no longer holds for weak spherical rank (upper or lower). The
main purpose of this section is to verify that claim.

\Small
\begin{table}[!ht]
   \begin{center}
    \begin{tabular}{| c | c | c |}\hline
    \multicolumn{3}{|c|}{\bf Compact manifolds} \\ \hline
     & $\forall \gamma$, there exists a& $\forall \gamma$, there exists a\\
    \textbf{Curvature} & Jacobi field $J$ s.t. & parallel vector field $E$
   s.t. \\
    \textbf{bound}& $\sect(J,\gamma')$ is extremal & $\sect(E,\gamma')$ is extremal. \\
    \hline \hline

     & &  the universal cover of $M$ is symmetric \\

    $\sect \leq 0$ & {\Large \textbf{?}} & or isometric to a product;\\
     & & {\Small cf.\ \cite{Ballmann}, \cite{BS}}.  \\
    \hline \hline

     & $M$ is isometric to & \\
    $\sect \leq -1$ & a locally symmetric space, &
    \multicolumn{1}{|c|}{$\boldsymbol{\Rightarrow}$} \\
     & {\Small cf.\ \cite{Hamenstadt}.} & \\
    \hline \hline

     & &\\
    $\sect \leq 1$ & non-symmetric examples exist. & $M$ is locally
    isometric to a CROSS.\\
     & & [ibid.]\\
    \hline \hline

     & & \\
    $\sect \geq 1$ & non-symmetric examples exist. & \multicolumn{1}{|c|}{\Large
    \textbf{?}} \\
     & & \\
    \hline \hline

    &  &{there are simply connected, irreducible}\\
    $\sect \geq 0$ & $\boldsymbol{\Leftarrow}$ & {examples which are
      not homeomorphic} \\

    & & to symmetric spaces, {\Small cf.\ \cite{Heintze},\cite{SS}.}\\
    \hline \hline

     & & \\
    $\sect \geq -1$ & {\Large \textbf{?}} &
    \multicolumn{1}{|c|}{{\Large \textbf{?}}} \\
     & & \\ \hline

    \end{tabular}
    \smallskip
   \end{center}
   \caption{Rank rigidity for various curvature bounds.}
   \label{Ta:first}

\end{table}
\normalsize

\subsection{The Berger spheres}
We present here, briefly, the construction of the so called Berger
spheres. This is the scaling of the round metric on a sphere; we will
specifically look at $\sph^3$. The Berger spheres are important
examples and originally were constructed by M.\ Berger in
\cite{Berger} to show that in odd dimensions, Klingenberg's
 injectivity radius
estimate  
 fails if the pinching is below
$\frac{1}{9}$.

One may regard the round 3-sphere as the unit sphere
in the quaternions  $\mathbf{H}$.
The Lie algebra is spanned  $i, j$ and $k$.
These vectors are orthonormal with respect
to a induced bi-invariant metric on $\sph^3$ of constant curvature $1$.

The Berger metric is obtained upon scaling the fibers
of the Hopf  fibration $\sph^3\rightarrow \sph^3/\sph^1=\sph^2$
where $\sph^1$ is the image of
the 1-parameter group $\exp(t i)$.
 More precisely, consider a family of left invariant metrics $g_{\eta}$
 on
$\sph^3$ which are defined by $g_{\eta}(i,j)=g_{\eta}(i,k)=g_{\eta}(k,j)=0$,
$\|j\|_{g_\eta}=\|k\|_{g_\eta}=1$ and $\|i\|_{g_\eta}=\eta$.

It is then routine to check that $\sect(\frac{1}{\eta}X_1, s\,X_2 +
\sqrt{1-s^2}\, X_3) = \eta^2$, $\sect(X_2, X_3) = 4-3\eta^2$, and
$\eta^2$ and $4-3\eta^2$ are minimum and maximum of the sectional
curvature; the Hopf fiber has length $2\pi \eta$
(for the calculation, see for instance \cite{CE}, Example 3.35).
 Note that if
$\eta>1$, then the range of curvatures is $[4-3\eta^2, \eta^2]$.

In order to find a non-symmetric example with weak lower spherical
rank take $\eta<1$ i.e., shrink the Hopf fiber and then normalize the
metric to make the lower bound 1. If $\gamma$ is a vertical geodesic
then all planes containing $\gamma'(t)$ have curvature 1.  If $\gamma$
is not vertical then the Killing field corresponding to the Hopf field
$i$ induces a Jacobi field $J$ along $\gamma$ with
$\sec(J,\gamma')=1$. Notice that $J$ is not necessarily normal but one
may replace $J$ by its normal part.

For weak upper spherical rank, take $\eta >1$ i.e., enlarge the
Hopf fiber (we may choose $\eta < \frac{2}{\sqrt{3}}$ to ensure
positive curvature). Normalize again to make the upper curvature bound
1 and as before it follows that the weak upper
spherical rank is 1.

Of course the Berger spheres in higher dimensions
also have positive weak upper or lower spherical rank.
Thus there are non-symmetric examples in all odd dimensions above 2.

It remains unclear whether the assumption on weak spherical rank
implies that the manifold is topologically a symmetric space. We leave
that as a question for further study.

\bibliographystyle{alpha}

\end{document}